\documentclass[12pt,a4paper]{amsart}
\usepackage[top=35mm, bottom=30mm, left=30mm, right=30mm]{geometry}
\usepackage{mathptmx}
\usepackage{mathrsfs}
\usepackage{verbatim}
\usepackage{url}
\usepackage[all]{xy}
\usepackage{xcolor}
\usepackage[colorlinks=true,citecolor=blue]{hyperref}
\usepackage{amsmath,amsthm}
\usepackage{amssymb}
\usepackage{enumerate}
\usepackage{graphicx}
\usepackage{tikz}
\usetikzlibrary{arrows}
\usetikzlibrary{patterns}

\newtheorem{thm}{Theorem}[section]

\newtheorem{prop}[thm]{Proposition}

\theoremstyle{definition}

\newtheorem{exa}[thm]{Example}

\newtheorem{conj}[thm]{Conjecture}

\numberwithin{equation}{section}
\begin{document}

%%%%% To ease editing, for IMPAN journals add:

\baselineskip=20pt

%%%%%%%%%%%%%%%%

\title{A note on the bounded orbit conjecture}

\author[E.~Shi]{Enhui Shi}
\address{School of Mathematics and Sciences, Soochow University, Suzhou, Jiangsu 215006, China}
\email{ehshi@suda.edu.cn}

%\author[H.~Xu]{Hui Xu*}
%\thanks{*Corresponding author}
%\address{ Department of mathematics, Shanghai normal university, Shanghai 200234,
% China}
%\email{huixu@shnu.edu.cn}

\author[Z.~Yu]{Ziqi YU}
\address{School of Mathematics and Sciences, Soochow University, Suzhou, Jiangsu 215006, China}
\email{20204207013@stu.suda.edu.cn}

\begin{abstract}
If $f:\mathbb{R}^2\rightarrow \mathbb{R}^2$ is an orientation reversing fixed point free homeomorphism on the plane $\mathbb{R}^2$ with no unbounded orbit, then $f$ has infinitely many periodic orbits.
\end{abstract}

\subjclass[2010]{}

\maketitle

\pagestyle{myheadings} \markboth{E. Shi and Z. Yu}{A note on the bounded orbit conjecture}

\section{Introduction}

In 1912, Brouwer \cite{Bro12b} proved his famous translation theorem as follows:

 \begin{thm}(Brouwer translation theorem)\label{brouwer}
 	If $f:\mathbb{R}^2\rightarrow \mathbb{R}^2$ is an orientation preserving homeomorphism and has no fixed point, then $f$ is a translation, that is each orbit of $f$ is unbounded.
 \end{thm}
\medskip

The following example shows that Brouwer translation theorem does not hold for orientation
reversing homeomorphisms (see e.g. \cite{Bo81}).

\begin{exa}\label{contex-Brouwer}
Let $h:\mathbb R^2\rightarrow \mathbb R^2, (x, y)\mapsto (-x, y)$ be the reflection
across $y$-axis. Let $g:\mathbb R^2\rightarrow \mathbb R^2$ be defined by $g(x, y)=(x, y-|x|+1)$
for $|x|<1$ and $g(x, y)=(x, y)$ for $|x|\geq 1$.
Then $g\circ h$ is an orientation reversing homeomorphism and has no fixed point, but each point $(x, y)$ with
$|x|\geq 1$ is a $2$-periodic point; specially, the orbit of $(x, y)$ is bounded.
\end{exa}

Noting that the homeomorphism $g\circ h$ in Example \ref{contex-Brouwer} has many unbounded orbits, the following conjecture is
natural.

\begin{conj}[Bounded orbit conjecture]
If $f:\mathbb R^2\rightarrow \mathbb R^2$ is a homeomorphism and each orbit of $f$  is bounded, then $f$
has a fixed point.
\end{conj}

 Boyles  solved this  conjecture in negative by constructing
a counterexample \cite{Bo81}. Recently, Mai etc. gave a more comprehensible and simpler counterexample \cite{MSYZ}.
We see that there are infinitely many periodic points in each of these counterexamples. The aim of the paper is to
show that this is always the case.

 \begin{thm}\label{main}
If $f:\mathbb{R}^2\rightarrow \mathbb{R}^2$ is an orientation reversing fixed point free homeomorphism  which has no unbounded orbit, then $f$ must have infinitely many periodic orbits.
\end{thm}

\section{Preliminaries}
In this section, we will recall some notions, notations, and elementary facts in the theories of
dynamical system and fixed point index theory.

\subsection{Some notions in dynamical system}

By a {\it dynamical system} we mean a pair $(X, f)$, where $X$ is a metric space
and $f:X\rightarrow X$ is a homeomorphism. For $x\in X$, the {\it orbit} of $x$ is
the set $O(x, f)\equiv\{f^i(x): i\in\mathbb Z\}$. If there is some $n> 0$ such that
$f^n(x)=x$, then $x$ is called a {\it periodic point} of $f$ and the minimal such
$n$ is called the {\it period} of $x$. A periodic point $x$ of period $1$ is called
a {\it fixed point}, that is $f(x)=x$. If there is a sequence of positive integers
$n_1<n_2<\cdots$ such that $f^{n_i}(x)\rightarrow x$, then we call $x$ a {\it recurrent point}.  Let ${\rm Fix}(f)$, $P(f)$ and $R(f)$ denote the sets of fixed points, period points and recurrent points respectively. A subset $A\subset X$ is called {\it invariant}  if $f(A)=A$. By an argument of Zorn's lemma, it is easy to see that every compact invariant set contains a recurrent point.
\medskip

Let $(X, f)$ be a dynamical system and $A$ be an invariant subset of $X$. For x$\in X$, if $$\lim_{n\rightarrow +\infty}d(f^n(x),A)\rightarrow 0,$$ then $x$ is said to be {\it positively asymptic to $A$}; if $$\lim_{n\rightarrow -\infty}d(f^n(x),A)\rightarrow 0,$$ then $x$ is said to be {\it negatively asymptic to $A$}; if $x$ is either positively asymptic to $A$ or negatively asymptic to $A$, then $x$ is said to be {\it asymptic to $A$}.
\medskip

Then following proposition can be found in \cite{GH}.

\begin{prop}\label{asy}
	Let $X$ be a compact metric space and $f:X\rightarrow X$ be a homeomorphism. If $A$ is a closed non-empty invariant set of $f$ and there exists a neighborhood $U$ of $A$ such that $x\in U-A$ implies $O(x,f)\not\subset U$, then there exists $y\in X-A$ such that $y$ is asymptotic to $A$.
	\end{prop}

\subsection{Fixed point index and Lefschetz number}

Let $W$ and $W'$ be two neighborhoods of 0 in the plane $(\mathbb{R}^2,\Vert\hspace{0.2em}\Vert)$ and $W$ contains the ball $B_r=\{(x,y)\in\mathbb{R}^2: \Vert (x,y)\Vert \leq r\}$ for some $r>0$. If $f: W\rightarrow W'$ is a homeomorphism having $0$ as the unique fixed point, then the index $i(0, f)$ is the degree of the map

$$\mu_0: S^1\rightarrow S^1, \quad s\mapsto \frac{f(rs)-rs}{\Vert f(rs)-rs \Vert}$$

where $S^1=\{s=(x,y)\in\mathbb{R}^2: \Vert s\Vert =1\}$. It is also the degree of the map

$$\mu: S^1\rightarrow S^1, \quad s\mapsto \frac{f(\gamma(s))-\gamma(s)}{\Vert f(\gamma(s))-\gamma(s) \Vert}$$

for any map $\gamma: S^1\rightarrow W-\{0\}$ homotopic to $\gamma_0 : S^1\rightarrow W-\{0\}$, $s\mapsto rs$.
\medskip

The fixed point index satisfies the following commutativity (see \cite{Dold}).

\begin{prop}\label{comm}
	Let $U$ and $V$ be two open sets in the plane $\mathbb{R}^2$. If maps $f: U\rightarrow \mathbb R^2$ and $g: V\rightarrow \mathbb R^2$
satisfy ${\rm Fix}(g\circ f)=\{a\}$ and ${\rm Fix}(f\circ g)=\{a'\}$ for some $a\in U$ and $a'\in V$, then $i(a,g\circ f)=i(a',f\circ g)$.
	\end{prop}

The following proposition is implied by the Lefschetz fixed-point formula (see \cite{Sp}).

\begin{prop}\label{lef}
	If $f: S^2\rightarrow S^2$ is an orientation preserving homeomorphism of the 2-sphere with only finitely many fixed points, then
	$\sum_{x\in Fix(f)}i(x,f)=2$.
	\end{prop}

The following result is due to Le Calvez \cite{Cal}.

\begin{prop}\cite[Proposition 1]{Cal} \label{rec}
Let $f: W\rightarrow W'$ be an orientation-preserving homeomorphism between two simply connected neighborhoods of $0\in \mathbb{R}^2$ such that ${\rm Fix}(f)=\{0\}$. Let $W''$ be the connected component of 0 in $W\cap f^{-1}(W')$. If there is a domain $V\subset W''$ such that $f(V)\cap V=\emptyset$, a point $z\in V$ and an integer $q\geq 2$ such that $f^i(z)\in W''$ for every $i\in\{1,\cdots,q-1\}$ and $f^q(z)\in V$, then $i(0,f)=1$.
\end{prop}

 \section{Proof of the main theorem \ref{main}}

\begin{proof}
	To the contrary, assume $f:\mathbb{R}^2\rightarrow \mathbb{R}^2$ is an orientation reversing fixed point free homeomorphism on the plane $\mathbb{R}^2$ which has no unbounded orbit and has only finitely many periodic orbits. Consider the orientation preserving homeomorphism $f^2:\mathbb{R}^2\rightarrow \mathbb{R}^2$. According to  Theorem \ref{brouwer}, ${\rm Fix}(f^2)\neq\emptyset$. Let $S^2=\mathbb R^2\cup \{\infty\}$ be the one point compactification of $\mathbb R^2$. Then we extend $f$ to a homeomorphism $\hat f$ on the sphere $S^2$ by letting $\hat f(\infty)=\infty$. Let $U$ and $V$ be simply connected neighborhoods of $\infty$ with $\infty\in \overline{U}\subset V$. Since $f$ has only finitely many periodic orbits, we can assume $\infty$ is the unique fixed point of $\hat f^2$ in $V$. As each orbit
of  $f^2$ is bounded, by Proposition \ref{asy}, there exists a point $p$ such that $\overline{O(p,f^2)}\subset\overline{U}\subset V$. Thus there is an $f^2$-recurrent point $x\in V$ which is not a fixed point. By Proposition \ref{rec} , the fixed point index $i(\infty, \hat f^2)$ is equal to $1$. For any fixed point $y$ of $f^2$, according to Proposition \ref{comm}, we have $i(y,f^2)=i(f(y),f^2)$. Thus the sum of all fixed point index of $\hat f^2: S^2\rightarrow S^2$ is an odd number. But by Proposition \ref{lef}, it must be $2$. This is a contradiction.
	\end{proof}

\subsection*{Acknowledgements} We would like to thank professor Jiehua Mai for helpful comments.

\end{document}